\documentclass[letterpaper, 10 pt, conference]{ieeeconf} 
\IEEEoverridecommandlockouts
\usepackage{float}
\usepackage{cite}
\usepackage{amsmath,amssymb,amsfonts}
\usepackage{algorithmic}
\usepackage{algorithm}
\usepackage{graphicx}
\usepackage{textcomp}
\usepackage{xcolor}
\usepackage{hyperref}
\usepackage{booktabs}
\hypersetup{
    colorlinks=true,
    linkcolor=black,
    filecolor=black,      
    urlcolor=blue,
    bookmarks=false,
    }
\def\BibTeX{{\rm B\kern-.05em{\sc i\kern-.025em b}\kern-.08em
    T\kern-.1667em\lower.7ex\hbox{E}\kern-.125emX}}
\usepackage[left=1.9cm, right=1.9cm, top=2.54cm]{geometry}
\title{\LARGE \bf
    Vehicle Routing for the Last-Mile Logistics Problem 
}

\author{Harshavardhan Desai, Ben Remer, Behdad Chalaki, Liuhui Zhao, \\Andreas A. Malikopoulos, Jackeline Rios-Torres 
\thanks{This research was supported in part by the Laboratory Directed Research and Development Program of the Oak Ridge National Laboratory, Oak Ridge, TN 37831 USA, managed by UT-Battelle, LLC, for the US Department of Energy (DOE), and in part by the Delaware Institute (DEI).}
    \thanks{Desai, Remer, Chalaki, Zhao and Malikopoulos are with the Department of Mechanical Engineering, University of Delaware, Newark, DE 19716 USA(e-mail:\texttt{harshd@udel.edu};\texttt{bremer@udel.edu}; \texttt{bchalaki@udel.edu};\texttt{lhzhao@udel.edu}; \texttt{andreas@udel.edu}}
    \thanks{Rios-Torres is with the Energy and Transportation Science Division, Oak Ridge National Laboratory, Oak Ridge, TN 37831 USA (e-mail: \texttt{riostorresj@ornl.gov)}}
    \thanks{This manuscript has been co-authored by UT-Battelle, LLC, under contract DE-AC05-00OR22725 with the US Department of Energy (DOE). The US government retains and the publisher, by accepting the article for publication, acknowledges that the US government retains a nonexclusive, paid-up, irrevocable, worldwide license to publish or reproduce the published form of this manuscript, or allow others to do so, for US government purposes. DOE will provide public access to these results of federally sponsored research in accordance with the DOE Public Access Plan (http://energy.gov/downloads/doepublic-access-plan).}%
}

\begin{document}

    \maketitle
    \thispagestyle{empty}
    \pagestyle{empty}

    \indent
    \begin{abstract}
        Energy consumption is the major contributor associated with large and growing transportation cost in logistics. Optimal vehicle routing approaches can provide solutions to reduce their operating costs and address implications on energy. This paper outlines a solution to the single-depot capacitated vehicle routing problem with the objective of minimizing daily operation cost with a homogeneous fleet of delivery vehicles. The problem is solved using Simulated Annealing, to provide optimal routes for the vehicles traveling between the depot and destinations. Simulation results demonstrate that the proposed approach is effective to recommend an optimal route and reduce operation cost. Supplementary information and video of our proposed approach can be found at: \url{https://sites.google.com/view/ud-ids-lab/last-mile}
    \end{abstract}

\section{Introduction}
    \subsection{Motivation}
    In a rapidly urbanizing world, fundamental transformations need to be made with respect to how transportation is used. We are currently witnessing an increasing integration of our energy and transportation networks which, coupled with the human interactions, is giving rise to a new level of complexity in transportation\cite{malikopoulos2015centralized}. With increasingly complex transportation systems\cite{malikopoulos2015duality}, new control approaches\cite{malikopoulos2015multiobjective},\cite{malikopoulos2018average} are needed to optimize the impact on system behavior of the interaction between vehicles at different applications.
    
    With the meteoric rise of the e-commerce industry, last-mile delivery, especially parcel delivery has attracted considerable attention\cite{hu2016us}. Last-mile logistics refers to the last portion of a supply chain involving the transportation of people or goods from the last transportation hub to the final destination. On the other hand, due to overcrowding and congestion in many metropolitan cities, more and more people are looking for ride-sharing as a viable form of transport for their daily commute\cite{burgstaller2017rethinking}. In this highly competitive environment, third party logistics and last-mile delivery firms must not only be able to meet ever-increasing fulfillment deadlines but do so as efficiently as possible. Some major hurdles in maximizing profits are ensuring the utilization of the vehicle inventory to its maximum potential and cost-effective routing of vehicles owing to growing emission constraints coupled with the strive to reduce fuel consumption.

    In this paper, we present a last-mile logistics system which combines two separate operations: (1) the transport of passengers using ride-sharing and (2) parcel delivery, using a single integrated fleet of vehicles, in order to minimize the total operational cost.

    \subsection{Literature Review}
    The vehicle routing problem (VRP) is a variation of the extensively investigated traveling salesman problem. There has been a significant amount of work done in the area of eco-VRP. Among different approaches that have been reported in the literature, time-dependent VRP (TD-VRP) is based on the notion that the travel time between any pair of points, such as customers and depots, depends on the distance between the points, as well as on the time of day (e.g., rush hours). The feature of fluctuating travel duration enables VRP to account for the actual conditions such as urban congestion, where the traveling speed is not constant due to variation in traffic density. Therefore, TD-VRP is a relevant and useful model to reveal recurring traffic congestion and to explore approaches to avoid it. In the model described by Qian and Eglese\cite{qian2016fuel}, the speed of the traffic on the underlying road network is time-dependent, and the path used by a vehicle between a pair of customers is the decision variable. The authors proposed a Tabu based algorithm to solve the problem and concluded that allowing a specified waiting time at customer nodes, and vehicles can avoid being caught in congestion, thus leading to overall fuel consumption reduction. Yao et al.\cite{yao2015vehicle} further explored the TD-VRP by introducing the concept of alternative stop points assuming that a delivery vehicle could temporarily stop at the opposite side of the client and then the deliveryman walks across the road to serve the client. By enabling alternative stops, detouring could be avoided in vehicle routing. As a result, the study showed that vehicle miles traveled as well as total fuel consumption were reduced in the network. Huang et al.\cite{huang2017efficient} presented an approach with the objective of finding an optimal routing solution such that vehicle arrival times at nodes meet the deadlines specified by the clients.
    
    Bent and Hentenryck\cite{bent2004scenario} considered the partially dynamic VRP with time windows. The goal was to serve as many customers as possible, given a fixed number of vehicles. For dynamic customers, stochastic information was assumed to be available. To tackle this dynamic stochastic vehicle routing problem, the authors proposed a multiple scenario approach, which involves continuously generating and solving scenarios with different static and dynamic requests, thereby generating an optimal routing plan.
    
    In the era of shared mobility, the ride-sharing problem has been widely studied. Furuhata et al.\cite{furuhata2013ridesharing} reviewed that key aspects of existing ride-sharing systems (e.g., the design of attractive and price effective mechanisms) along with some of their key challenges such as ride-arrangement customer preferences, multi-modal rides, and building of trust among unknown travelers. Agatz et al.\cite{agatz2012optimization} surveyed the characteristics, objectives, and optimization challenges of different classes of operations research models related to ride-sharing. The cost of the ride-share trip should be proportionately divided among the participants, roughly proportional to vehicle-miles traveled. Agatz et al.\cite{agatz2011dynamic} proposed a framework allocating the costs of a joint trip proportional to the distances traveled if separate trips were taken. Kleiner et al.\cite{kleiner2011mechanism} proposed an auction-based mechanism to determine the driver’s compensation.
    
    Wang et al.\cite{wang2016towards} presented a model which utilizes crowd-workers for last-mile parcel delivery. The concept was such that each parcel will be sent to a pick-own-parcel station nearest to its consumer’s address, then assigned to the crowd-workers via a mobile app. Each worker was associated with a travel pattern (e.g., the driver’s daily commute route) and the reward was proportional to the distance of the detour taken to make a delivery. The objective was to assign the parcels to the most convenient workers to minimize the total reward paid by the company. Compared to traditional delivery methods, the crowd-sourcing approach resulted in a higher level of parallelism in job execution since the fleet size of delivery vehicles is much larger and each worker handles only a small number of parcels. Communication between workers and customers is also more effective. Due to the elimination of vehicles specifically used for last-mile delivery, this approach helped reduce operations costs as well as carbon emissions, since the number of vehicles in the network were reduced.

    \subsection{Contribution of this Paper}
    Both ride-sharing and last-mile delivery are topics that have been extensively studied with either static or dynamic settings. However, the concept of shared logistics has not been fully investigated in the literature. The main objective of this paper is to show how a fleet of vehicles can be dispatched most efficiently to maximize the daily profit of a company which aims at fulfilling the last/first-mile connection of passengers and simultaneous parcel deliveries from a major transport hub. By doing so, we also show how much more efficient this concept can be over the baseline scenario (using two dedicated fleets). The concept is then also replicated on a scaled testbed.
    \subsection{Organization of the Paper}
    The structure of the paper is organized as follows. In Section II, we formulate the problem. In Section III, we provide the solution and simulation results. Finally, in Section IV, we draw concluding remarks.
    
\begin{figure}
    \centering
    \includegraphics[width=0.45\textwidth]{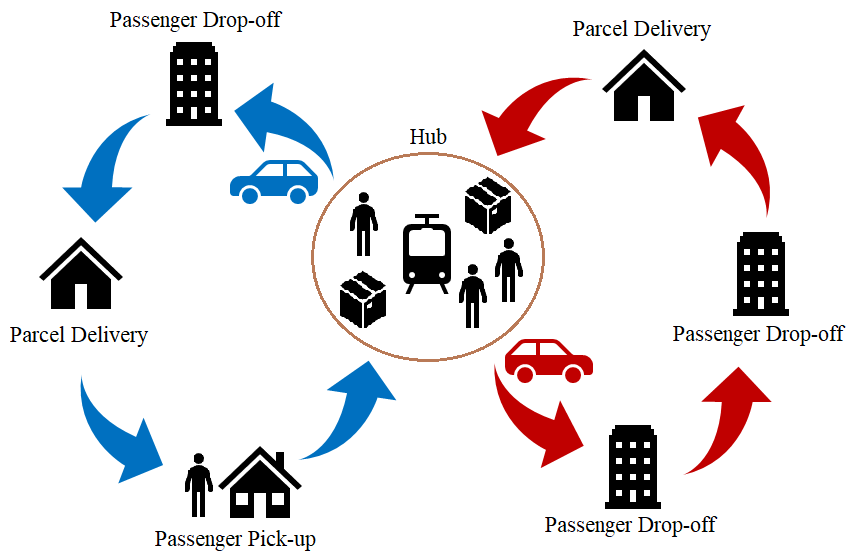}
    \caption{Last-mile city logistics problem.}
    \label{fig:Figure 1}
\end{figure}   

\section{Problem Formulation}
Consider a major transit hub in a city, which serves as a consolidation center for parcels to be delivered to residents in the neighborhoods. We consider a first/last-mile transportation service for passengers stationed around the hub as well as for the delivery of parcels in the same area. Electric autonomous vehicles owned by the taxi/ride-hailing company are used for simultaneous transport and delivery. A representation of the problem is shown in Fig. \ref{fig:Figure 1}.
    \subsection{Problem Description}
    \textbf{Transportation of Passengers:} All vehicles are located at the transport hub at the start of the operation time where parking and charging are free. Similarly, all vehicles return to the hub after the completion of their route. This guarantees that the batteries of the vehicle fleet can be charged and any required maintenance can be performed. The pick-up or drop-off locations of the passengers are known beforehand. With this information, optimal routes for a fleet of vehicles with a fixed capacity, which serve all the passengers, can be generated.

\begin{figure}
    \centering
    \includegraphics[width=0.3\textwidth]{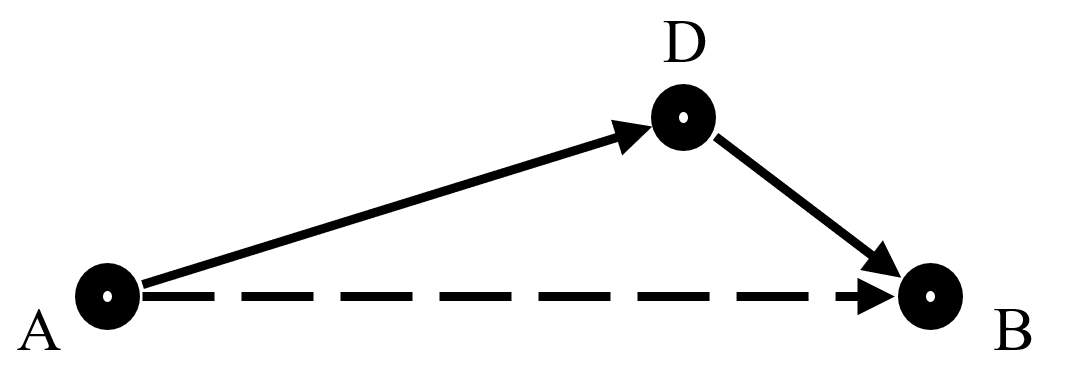}
    \caption{Example of a delivery path.}
    \label{fig:Figure 2}
\end{figure}

    \textbf{Parcel Delivery:} With our proposed model, the service company only needs to handle the delivery of parcels to the hubs nearest to the consumer’s address after which, the parcels at the hub will be assigned to a vehicle and eventually reach the consumers. A parcel will be assigned to that vehicle which will incur the least cost deviating from its assigned route to make the delivery. For example, consider a vehicle traveling from node A to node B, as shown in the dashed line in Fig. \ref{fig:Figure 2}. For simplicity, assume that the cost, $r$, of deviating from the path and making a delivery at node D on the way, is proportional to the additional travel distance
    
    \begin{equation} \label{eq:traveldistance}
    r=d(A,D)+d(D,B)-d(A,B),
    \end{equation}
    where $d(\cdot,\cdot)$ is the Euclidean distance between two locations, $d(A,D)+d(D,B)$ is the distance incurred for the delivery, and $d(A,B)$ is the travel distance if the vehicle does not take the task. 
    
    \textbf{Routing:} Consider a network with $N+1$, $N\in\mathbb{N},$ nodes. Let $\mathcal{N}=\{1,\ldots,N\}$ be the customer nodes and \textbf{\textit{0}} be the node that represents the hub.

    \textit{Routes:} A route indicates a path for a vehicle that starts at the hub, serves a set number of customers, and then returns to the hub. Mathematically, it is a sequence of nodes, $\{0,s_1,\dots,s_n,0\}$, where $s_i$ is a node, $1\leq i\leq N$, and all $s_i$ are distinct. The travel cost of a route $r$ is denoted by $c(r)$ and is the summed cost of visiting all of its customers, i.e., $c(r)=c_{0s_1}+c_{s_1 s_2}+\dots+c_{s_{n-1} s_n }+c_{s_n 0}$.
    
    \textit{Routing Plan:} A routing plan, is a set of routes ${r_1,\dots,r_m }$  serving each customer exactly once. The travel cost of a plan is denoted by $c(\alpha)=\Sigma_{i=1}^n c(r)$.
    
    In the modeling framework above, the following assumptions are imposed:

    \textbf{Assumption 1:} A vehicle can accommodate up to two passengers and one parcel, or one passenger and three parcels at a time.
    
    \textbf{Assumption 2:} The vehicles are fully charged overnight at the transport hub.
    
    These assumptions consider the capacity constraint for each vehicle and ensure that the trip is short enough to be completed on a single charge. They also ensure that the vehicles are fully charged at the start of the operational period.
    \subsection{Objective Function}
    The revenue earned from carrying passengers from the transport hub to their destinations is
    
    \begin{equation}\label{eq:revenuefromTH}
        P\cdot\sum_{i\in \mathcal{N} }D_{0i}d_{0i},
    \end{equation}
    where $P$ is the price rate per driving distance (\$/mile), $d_{0i}$ is the travel distance between the transport hub and service node $i$, $D_{0i}$ is the number of trips satisfied from the hub to node for all ${i\in \mathcal{N}}.$
    
    The revenue earned from carrying passengers from their origins to the transport hub is
    \begin{equation}\label{eq:revenuetoTH}
        P\cdot\sum_{i\in \mathcal{N} }D_{i0}d_{i0},
    \end{equation}
    where $d_{i0}$ is the travel distance between service node $i$ and the transport hub, $D_{i0}$  is the number of trips satisfied from node $i$ to the hub, for all ${i\in \mathcal{N}}.$
    
    The total vehicle running cost is
    \begin{equation}\label{eq:runningcosts}
        C_d\cdot\sum_{i,j\in \mathcal{N}, \; i\neq j}U_{ij}d_{ij},
    \end{equation}
    where $C_d$ is the vehicle running costs per mile (\$/mile), $U_{ij}$ is the number of vehicles travelling from node $i$ to node $j$.
    
    The total vehicle maintenance cost is
    \begin{equation}\label{eq:maintenance}
        C_mF,
    \end{equation}
    where $C_m$ is the maintenance cost per vehicle per day (\$/day), $F$ is the vehicle fleet size in the system.
    
    The total detour cost for delivering parcels is
    \begin{equation}\label{eq:detour}
        \sum_{i,j\in \mathcal{N}, \; i\neq j}x_{ij}d_{ij}C_p,
    \end{equation}
    where $x_{ij}$ is equal to $1$, if there is a vehicle travelling from node $i$ to node $j$, otherwise $0$, for all ${i,j\in \mathcal{N}}, d_{ij}$ is the travel distance between node $i$ and node $j$, for all ${i,j\in \mathcal{N}},C_p$ is the detour cost per mile (\$/mile).
    
    The following problem formulation seek to maximize the total profit, $G$, during a typical day of operations, by summing the revenues earned in \eqref{eq:revenuefromTH}, \eqref{eq:revenuetoTH} and the above costs  \eqref{eq:runningcosts}-\eqref{eq:detour}.
    \begin{equation}\label{eq:Objective}
    \begin{aligned}
        \max_{D_{0i}, D_{i0},U_{ij},x_{ij}}{G}=
        P\cdot\bigg(\sum_{i\in \mathcal{N} }D_{0i}d_{0i}+\sum_{i\in \mathcal{N} }D_{i0}d_{i0}\bigg)\\
        -C_d\cdot\sum_{i,j\in \mathcal{N}, \; i\neq j}U_{ij}d_{ij}-C_mF-\sum_{i,j\in \mathcal{N}, \; i\neq j}x_{ij}d_{ij}C_p
        \end{aligned}
    \end{equation}

    \subsection{Constraints}
    The following constraints describe conservation of flow in the network
    \begin{equation}\label{constraint:flow}
        \sum_{j\in \mathcal{N}, \; i\neq j}U_{ji} =\sum_{j\in \mathcal{N}, \; i\neq j}U_{ij},
    \end{equation}
     which states that the number of vehicles entering a node $i$ from any other node, must be equal to the number of vehicles exiting node $i$.
    
    The following constraint ensures that demand is satisfied
    \begin{equation}\label{constraint:demand}
        D_{0i}\leq Q_{0i}x_i, \quad \forall i \in \mathcal{N},
    \end{equation}
    where, $x_i$ is $1$ if the passenger request at node $i$ can be served, $0$ otherwise, for all ${i\in \mathcal{N}}$, $Q_i$ is the number of passenger requests from the transport hub to node $i$. The constraint (\ref{constraint:demand}) assures that the trips satisfied between the transport hub and service node $i$ must be lower than or equal to the passenger requests on the same origin-destination (OD). If node $i$ cannot be served ($x_i=0$), the satisfied demand must be zero.
    
    The following constraint assures that the trips satisfied between service node $i$ and the transport hub must be lower than or equal to the passengers’ requests on the same OD pair. If node $i$ cannot be served ($x_i=0$), the satisfied demand must be zero, namely
    \begin{equation}\label{constraint:odpair}
        D_{i0} \leq Q_{i0}x_i, \quad \forall{} i \in \mathcal{N},
    \end{equation}
    where $Q_i$  is the number of passenger requests from node $i$ to the transport hub.
    
    We also need to ensure that if no trip from node $i$ is satisfied, then that node is not selected, hence
    \begin{equation}\label{constraint:notrip}
        x_i\leq\sum D_{0i} +\sum D_{i0}, \quad \forall i \in \mathcal{N}.
    \end{equation}

    Next, we need to include a constraint that imposes the condition that the number of vehicles traveling between service node i and the transport hub must be greater than or equal to the number of passengers traveling on that OD route
    \begin{equation}\label{constraint:passengers}
        D_{i0} \leq U_{i0}, \quad \forall i \in \mathcal{N}.
    \end{equation}
    
    Finally, the next set of constraints define the domain for the decision variables,
    \begin{equation}
        D_{0i} \geq 0, \quad \forall i \in \mathcal{N},
    \end{equation}
    \begin{equation}
        D_{i0} \geq 0, \quad \forall i \in \mathcal{N},
    \end{equation}
    \begin{equation}
        U_{ij} \geq 0, \quad \forall i,j \in \mathcal{N},
    \end{equation}
    \begin{equation}
        x_i,x_{ij}\in(0,1), \quad \forall i \in \mathcal{N}.
    \end{equation}
\section{Solution Approach and Simulation Results}

\begin{table*}[!ht]
    \caption{Optimal Results for ten Instances}
    \label{tab:teninstance}
    \centering
 \begin{tabular}{cccccccc} 
        Instance&No. of &No. of &No. of& Total no. of &Baseline Cost&Cost with &Reduction in Cost\\
        &Vehicles&Passengers &Deliveries &Customers &(\$ USD)&Integrated Service (\$ USD)& (\%)\\
         \toprule
         \textbf{1}&3&6&2&8&461.68&246.77&46.55\\
         \textbf{2}&4&7&3&10&425.65&231.11&45.71\\
         \textbf{3}&5&10&4&14&656.31&290.36&55.76\\
         \textbf{4}&7&14&6&20&500.87&303.15&39.48\\
         \textbf{5}&9&17&8&25&592.08&326.07&44.93\\
         \textbf{6}&11&20&10&30&624.03&363.01&41.83\\
         \textbf{7}&14&27&13&40&692.06&466.04&32.66\\
         \textbf{8}&17&34&16&50&738.08&501.74&32.02\\
         \textbf{9}&21&40&20&60&743.07&521.61&29.80\\
         \textbf{10}&24&47&23&70&802.73&544.70&32.14\\
\midrule
    \end{tabular}
\end{table*}

\begin{algorithm}
    \caption{Pseudocode for Simulated Annealing}
    \label{alg:Solver}
    \begin{algorithmic}
        \STATE{\textbf{Input:} ProblemSize, iterations\_max, temp\_max\\}
        \STATE{\textbf{Output: }\textit{: S\_best}\\}
        \STATE{$S\_c\longleftarrow CreateInitialSolution(ProblemSize)$\\}
        \STATE{$S\_best\longleftarrow S\_c$\\}
        \FOR{$(i = 1$ To $iterations\_max)$}
        \STATE{$S\_i \longleftarrow CreateInitialSolution(S\_c)$\\}
        \STATE{$temp\_c \longleftarrow CalcTemperature(i, temp\_max)$\\}
        \IF{$(Cost(S\_i)\leq Cost(S\_c))$}
        \STATE{$S\_c \longleftarrow S\_i$\\}
        \IF{$(Cost(S\_i)\leq Cost(S\_best))$}
        \STATE{$S\_best \longleftarrow S\_i$\\}
        \ENDIF{\\}
        \IF{$(Exp((Cost(S\_c)– Cost(S\_i))/temp\_c)>Rand())$\\}
        \STATE{$S\_c \longleftarrow S\_i$}
        \ENDIF{\\}
        \ENDIF{\\}
        \ENDFOR{\\}
        \RETURN{$(S\_best)$}
    \end{algorithmic}
\end{algorithm}

The VRP is classified as an NP-hard problem. When the problem involves real-world data sets that are usually very large, it may become difficult to solve the problem within acceptable CPU times, if exact optimization methods such as direct tree search methods, dynamic programming, and integer linear programming are used. The vehicle routing problem comes under combinatorial problem. Hence, to get solutions in determining routes which are realistic and very close to the optimal solution, heuristics and meta-heuristics are used. Simulated Annealing (SA), Tabu Search, Genetic Algorithm and Ant Colony Optimization are some of the meta-heuristics that have been applied to the VRP.

SA is inspired by the process of annealing in metallurgy. In this natural process, a metal is heated and slowly cooled under controlled conditions to increase the size of the crystals in the metal and reduce its defects. The heat increases the energy of the atoms allowing them to move freely, and the slow cooling schedule allows a new low-energy configuration to be discovered and exploited. Similarly, each configuration of a solution in the search space represents different internal energy of the system. Heating the system results in a relaxation of the acceptance criteria of the samples taken from the search space. As the system is cooled, the acceptance criteria of samples are narrowed to focus on improving movements. Once the system has cooled, the configuration will represent a sample at or close to a global optimum.

The main idea of a SA algorithm is to occasionally accept degraded solutions in the hope of escaping the current local optimum. The information processing objective of the technique is to locate the minimum cost configuration in the search space. The algorithm’s plan of action is to probabilistically re-sample the problem space where the acceptance of new samples into the currently held sample is managed by a probabilistic function that becomes more discerning of the cost of samples it accepts over the execution time of the algorithm. The pseudocode for the SA algorithms used for the VRP is given in Algorithm 1, codes for the algorithm can be found online at \cite{yarpiz,lzane}.


\begin{figure}
    \centering
    \includegraphics[width=0.45\textwidth]{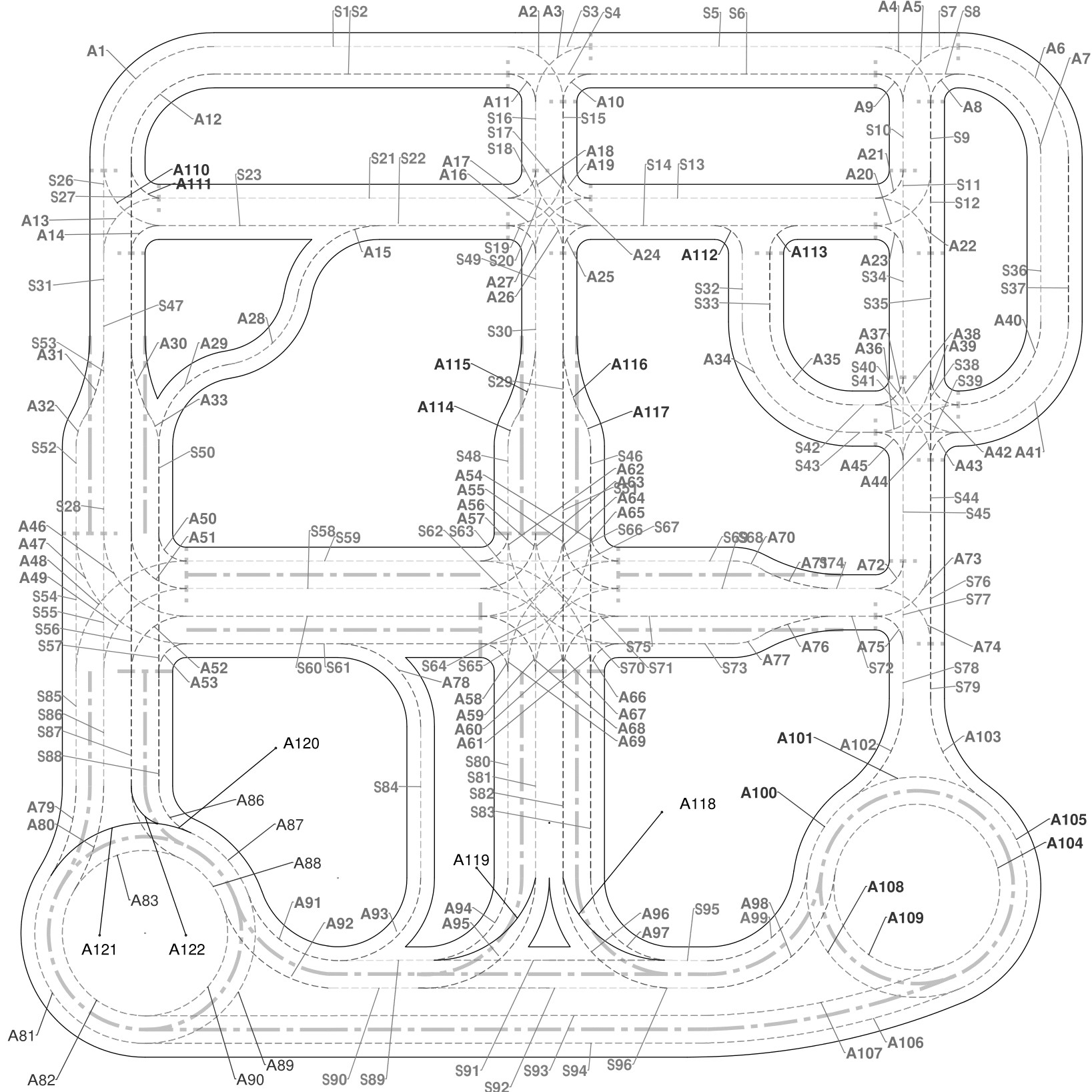}
    \caption{University of Delaware Scaled Smart City network.}
    \label{fig:Figure 3}
\end{figure}

To evaluate the performance of the proposed approach, ten computational experiments were carried out on a network spanning 0.25‬ sq.km. (shown in Fig. \ref{fig:Figure 3} the University of Delaware Scaled Smart City network) with a varying number of customer locations as well as varying fleet sizes. The number of customer locations vary from 8 to 70, out of which two-thirds of the locations are passenger locations, and the rest are parcel delivery locations. The capacity of a vehicle is limited to 5 units, where each passenger takes two units, and each parcel takes 1 unit. This is done to ensure that no more than two passengers are seated in the vehicle at a time. With one passenger, there will be room for three parcels, and with two passengers, there will be room for one parcel.

For the baseline scenario, we have considered a combination of two conventional last-mile transportation system, which consists of two separate fleets of vehicles, one for passenger transport (ride-sharing) and one for parcel delivery. The operational costs in this scenario would be the total costs of operating both fleets simultaneously. For simplicity, we consider the vehicle type and capacity of the vehicles of both fleets to be the same. Keeping the vehicle capacity in mind, in the baseline scenario, the passenger fleet size varies from 3 to 24 vehicles, and the delivery fleet size varies from 1 to 5 vehicles, proportionate to the number of customer locations, bring in the total vehicle count to 4 to 29 vehicles.

With the integrated service, the fleet size varies from 3 to 24 vehicles, proportionate to the number of customers (i.e., 8 to 70 customer locations). We use SA to find the cost of every possible route that can be taken for every instance until the minimum cost (and its corresponding vehicle routes) is obtained. The routing plan, which leads to the lowest cost, is taken as the optimal solution. The computation of the SA was executed on MATLAB running on an computer (Intel Core $i7-6700$ CPU @ $3.40$ GHz) running Windows 10.

\begin{figure}
    \centering
    \includegraphics[width=0.5\textwidth]{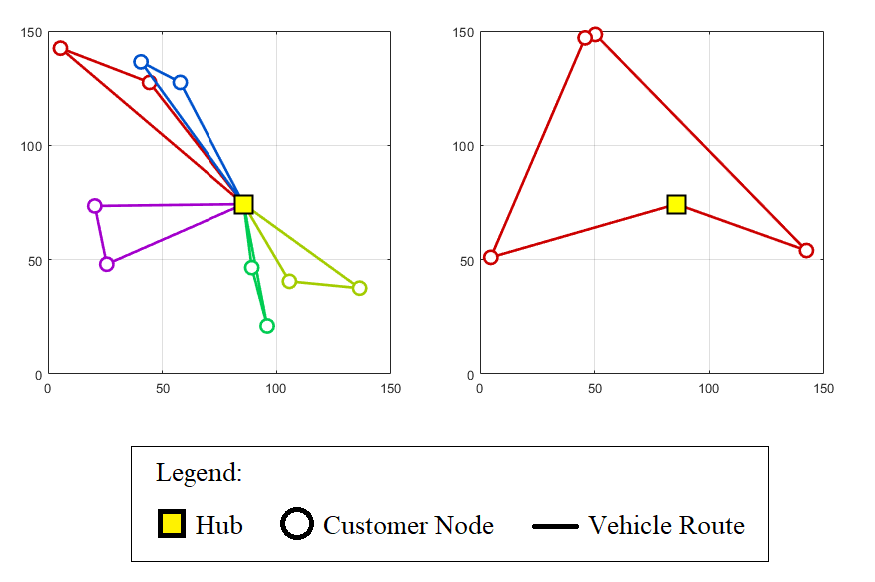}
    \caption{Baseline scenario for instance 3 (Left: Passenger fleet, Right: Delivery fleet).}
    \label{fig:Figure 4}
\end{figure}

\begin{figure}
    \centering
    \includegraphics[width=0.4\textwidth]{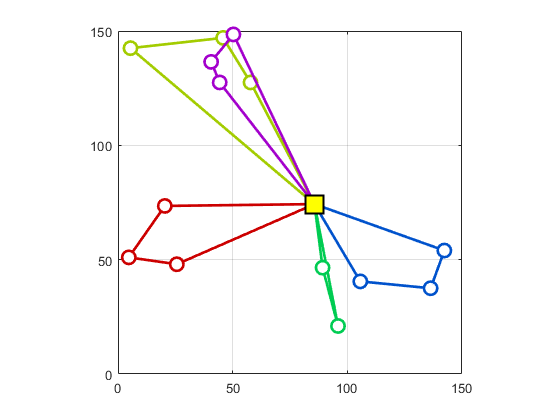}
    \caption{Optimal Routes for instance 3.}
    \label{fig:Figure 5}
\end{figure}

\begin{figure}
    \centering
    \includegraphics[width=0.5\textwidth]{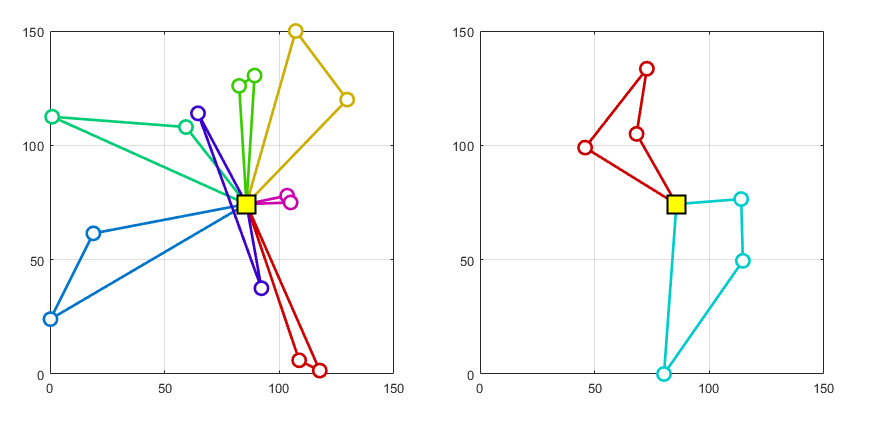}
    \caption{Baseline Scenario for instance 4 (Left: Passenger fleet, Right: Delivery fleet).}
    \label{fig:Figure 6}
\end{figure}

\begin{figure}
    \centering
    \includegraphics[width=0.4\textwidth]{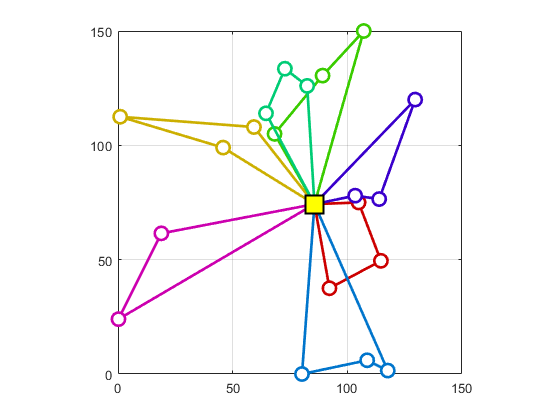}
    \caption{Optimal Routes for instance 4.}
    \label{fig:Figure 7}
\end{figure}

\begin{figure}
    \centering
    \includegraphics[width=0.4\textwidth]{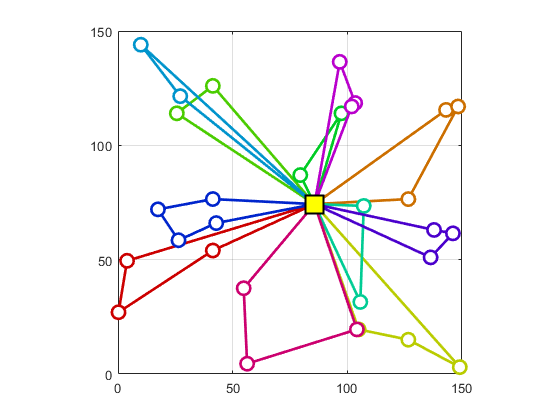}
    \caption{Optimal Routes for instance 6.}
    \label{fig:Figure 8}
\end{figure}

\begin{figure}
    \centering
    \includegraphics[width=0.4\textwidth]{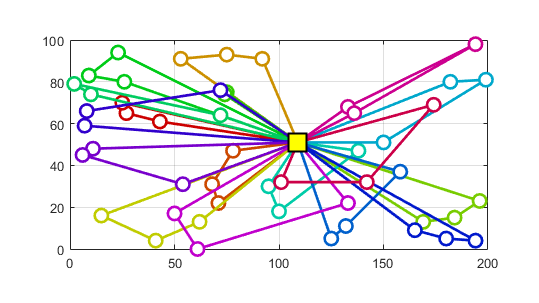}
    \caption{Optimal Routes for instance 8.}
    \label{fig:Figure 9}
\end{figure}

\begin{figure}
    \centering
    \includegraphics[width=0.4\textwidth]{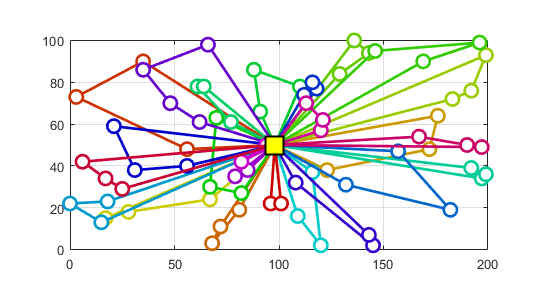}
    \caption{Optimal Routes for instance 10.}
    \label{fig:Figure 10}
\end{figure}

The experimental results for the ten instances are presented in Table I. It can be observed that the total cost of service is roughly proportional to the number of customer locations and fleet size. For example, the cost of serving 25 customer locations using nine vehicles (instance 5) is higher than the cost of serving 20 customer locations using seven vehicles (instance 4). But it is also observed that the cost of servicing ten customers using four vehicles (instance 2) is greater than the cost of servicing 14 customers using five vehicles (instance 3). This behavior can be due to the fact that in instance 2, one of the vehicles serves only one customer, indicating that the fleet is underutilized. The optimal routes for instances 3, 4, 6, 8, and 10 are shown in Fig. \ref{fig:Figure 4}-\ref{fig:Figure 10}, in which we see that in each route, each vehicle visits two passenger locations for drop-off or pick-up, as well as a parcel delivery location before returning to the hub.

We observe that compared to the baseline scenario, the integrated service reduces operational costs by a considerable amount. Cost reduction ranges from 29.8\% during high demand (instance 10), to 55.76\% during low demand (instance 1). Take instance 3 for example: Since the delivery locations are spread out, the sole delivery vehicle must travel a large amount to make all four deliveries, resulting in a high operating cost (\$386.09). The passenger fleet costs amounted to \$270.22 to cater to 10 passengers, bringing the total cost of using both fleets to \$656.31. The vehicle routes in instance 3 using the baseline scenario is shown in Fig. \ref{fig:Figure 4}. Comparatively, in the integrated service, when the passenger vehicles are used to make the deliveries, the cost of operating the delivery vehicle is eliminated and is offset by a marginal increase in the operating costs of the passenger fleet, since they must travel a bit more to make the nearby delivery. The operating cost for this instance was \$290.36, which is 55.76\% lower than the baseline scenario.

\section{Experimental Deployment}

\begin{table*}[!ht]
\caption{Experimental Results}
     \label{table:ExRes}
\centering
\begin{tabular}{cc|cc|cc}
     Instance&No. of CAVS&\multicolumn{4}{c}{Time Penalty Over Passenger Fleet(Min.)} \\
     &&Total (original)&Total (scaled)&Average (original)&Average (scaled)\\
\toprule     1&3&0.36&\bf{8.96}&0.12&\bf{2.99}\\
     2&4&0.37&\bf{9.21}&0.09&\bf{2.30}\\
    \hline{}
\end{tabular}
\end{table*}

\begin{figure}
    \centering
    \includegraphics[width=0.45\textwidth]{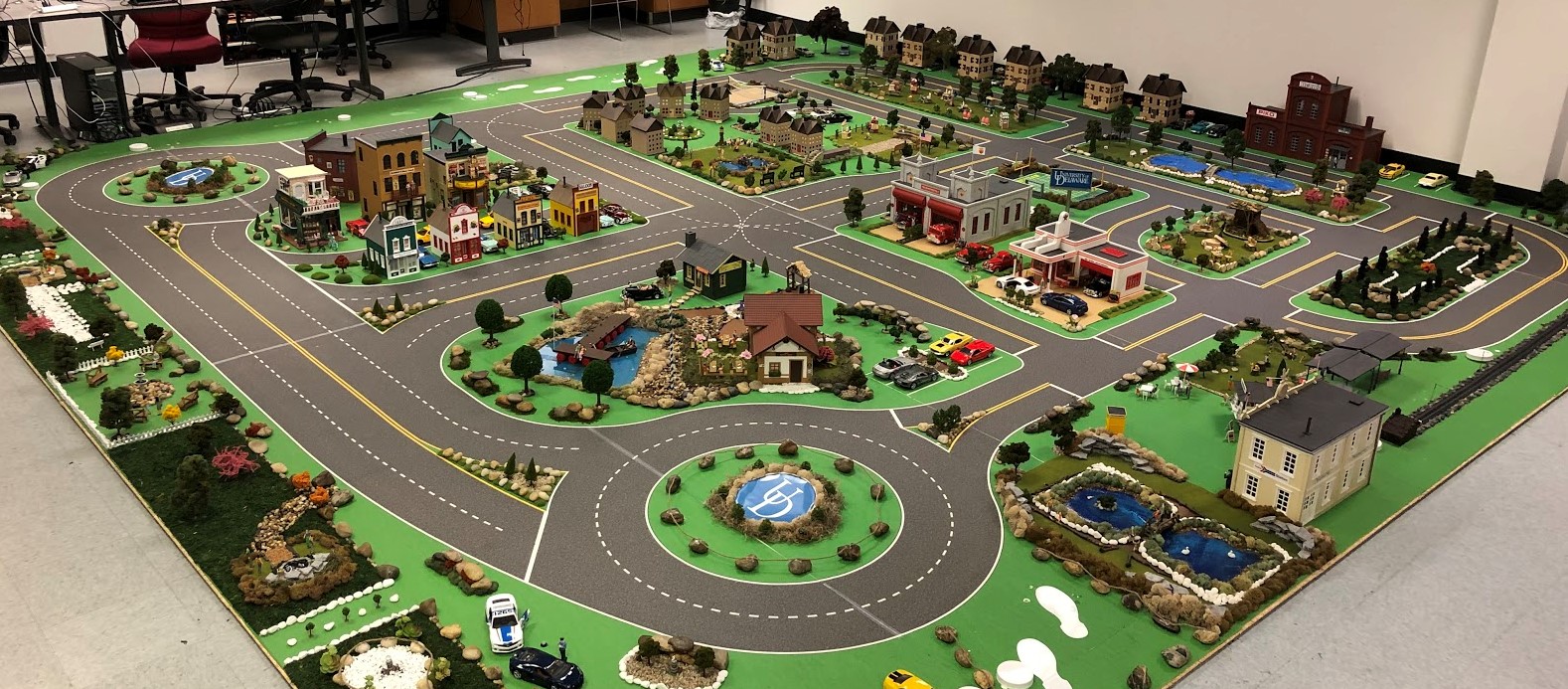}
    \caption{University of Delaware Scaled Smart City (UDSSC).}
    \label{fig:Figure 12}
\end{figure}

\begin{figure}
    \centering
    \includegraphics[width=0.45\textwidth]{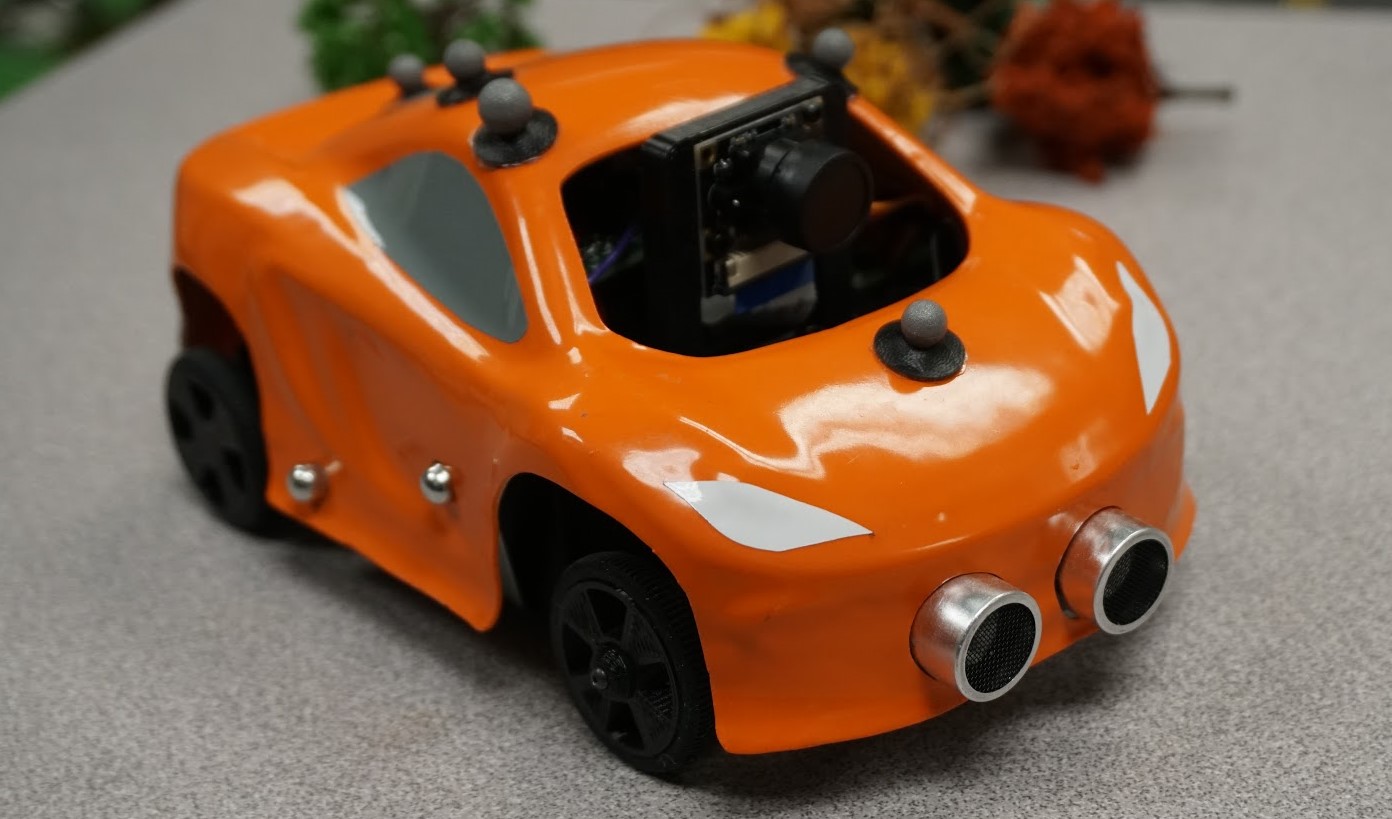}
    \caption{CAVs used for the experiments.}
    \label{fig:Figure 13}
\end{figure}

To further validate the concept, two experiments were carried out on the University of Delaware Scaled Smart City (UDSSC), a 1:25 scale testbed (Fig. \ref{fig:Figure 12}), designed to replicate real-world traffic scenarios and test cutting-edge control technologies in a safe and scaled environment. UDSSC is a fully integrated smart city, which can be used to validate the efficiency of control and learning algorithms and their applicability in hardware. It utilizes high-end computers, a VICON motion capture system, and scaled CAVs (Fig. \ref{fig:Figure 13}) to simulate a variety of control strategies with up to 35 scaled CAVs. Each CAV has a Raspberry Pi 3B with a $1.2$ GHz quad-core ARM processor and communicates with the \textit{mainframe} computer (Processor: Intel Core $i7-6950X$ CPU @ $3.00$ GHz x $20$, Memory: $125.8$ Gb). UDSSC has been used successfully for coordination of CAVs \cite{Malikopoulos2018b, beaver2019demonstration} and implementation of reinforcement learning policies \cite{jang2019simulation, chalaki2019zero}.

In this paper, the CAVs acted like the vehicles used for \textit{passenger ride-sharing} and \textit{parcel delivery}. The CAVs ran at a speed of 0.4 m/s, which when scaled to real-world conditions, would be around 17 mph, which would be a reasonable speed for a residential area.

Instances 1 and 2 from the simulations were chosen to be the basis for the experiments. The results from the simulation gave us the optimal routing plans, containing the optimal routes for each vehicle, for both the baseline scenarios as well as for the integrated service. The routes were then assigned to the vehicles, using an Intelligent Driver Model (IDM) to command the movement of every CAV on the scaled city. The IDM ensures that the CAVs follow the desired path at the desired speed, and takes preventive measures to avoid collisions. First, the baseline scenarios were demonstrated, where two separate fleets were used for deliveries and passenger. Then the integrated service was demonstrated, by using the passenger fleet to carry out the deliveries, as assigned by the simulation results. The implication of this is that the CAVs would be operational for a bit longer as compared to the baseline scenario of the passenger fleet, but in doing so, we eliminate the delivery fleet altogether. Videos of the experiments along with the corresponding simulations and the actual routes taken are posted on: \url{https://sites.google.com/view/ud-ids-lab/last-mile}

The results in Table II show the time penalty that the fleet for the integrated service suffers over the passenger fleet in the baseline scenario. According to the results, the total time that the CAVs were operational was around 0.36 minutes more, which when scaled to real-world conditions equates to around 9 minutes more, than the amount of time for the pure passenger fleet. On average, every CAV was operational for 2-4 minutes more. This penalty is acceptable since it allows for the delivery fleet to be completely eliminated.

\section{Conclusions}
The VRP is an interesting problem not only for distribution centers but also for shared mobility services. This problem has been addressed in the literature. However, combining the two areas of VRPs has not been fully investigated. In this paper, we proposed an efficient vehicle routing planning scheme with the objective of minimizing the daily cost of operations of a fleet of vehicles, which is used to ferry passengers to their final destination, and simultaneously deliver parcels to nearby destinations. The effectiveness of the approach was tested on a benchmarking network using a SA algorithm. The results indicated that the proposed approach is valid, reliable, and has good computing performance. It was shown that the integrated service can reduce operational costs by up to 56\%, compared to conventional last-mile transportation services, depending on customer demand. Furthermore, the experimental validation at the UDSSC testbed showed that the concept is practical and feasible in real world conditions.

In the computational model, only a few factors were considered, namely demand at passenger and delivery locations, vehicle capacity, and transportation costs. As it is known, fuel efficiency is highly sensitive to the driving cycle, this method can be further enhanced with dynamic traffic flow data so that fuel consumption can be reduced. Finally, this approach can be modified to accommodate a heterogeneous delivery fleet (i.e., vehicles with different cargo-carrying capacities as well as vehicles with different kinds of powertrains).

\bibliographystyle{IEEEtran}
\bibliography{main}

\end{document}